\documentclass[article]{amsart}
\usepackage{hyperref}
\numberwithin{equation}{section}
\newtheorem{theorem}{Theorem}[section]

\newtheorem{lemma}[theorem]{Lemma}

\begin{document}
\title[{Meromorphic functions share three values} ]{Meromorphic functions share three values with their difference operators }
\author[F. L\"{u} and W.R. L\"{u}]{Feng L\"{u} and Weiran L\"{u}}
\address{College of Science\\
China University of Petroleum\\
Qingdao, Shandong, 266580, P.R. China.} \email{lvfeng18@gmail.com}

\address{College of Science\\
China University of Petroleum\\
Qingdao, Shandong, 266580, P.R. China.} \email{luwr@upc.edu.cn}

\thanks{The research was supported by the Natural Science
Foundation of Shandong Province Youth Fund Project (ZR2012AQ021)}


\begin{abstract}
In the work, we focus on a conjecture due to Z.X. Chen and H.X. Yi \cite{CY} which is concerning the uniqueness problem of meromorphic
functions share three distinct values with their difference
operators. We prove that the conjecture is right for meromorphic function of finite order. Meanwhile, a
result of J. Zhang and L.W. Liao \cite{ZL} is generalized from
entire functions to meromorphic functions.
\end{abstract}

\maketitle

\section{Introduction and main result}
In Nevanlinna theory, to consider the relationship of two
meromorphic functions if they share several values CM or IM is an
important subtopic, such as the famous Nevanlinna's five and four
values theorems \cite{N}. In 1976, L.A. Rubel and C.C. Yang \cite{RY}
showed that if nonconstant entire function $f$ and its first
derivative $f'$ share two distinct values CM, then they are
identical. This result is extended by E. Mues and N. Steinmetz
\cite{MS} in 1979 from shared values CM to IM, by L.Z Yang \cite{Y}
in 1990 from first derivative to $k$-th derivative.\\

As the difference analogues of Nevanlinna¡¯s theory are being
investigated, J. Zhang and L.W. Liao \cite{ZL} considered the
difference analogues of the result of L.A. Rubel and C.C. Yang. They
replaced the derivative $f'$ by the difference operator $\Delta f=
f(z+1)-f(z)$ and obtained the following result.\\

{\bf Theorem A}. {\it Let $f$ be a transcendental entire function of
finite order and $a,~ b$ be two distinct constants. If $\Delta f (\not\equiv0)$
and $f$ share $a,~ b$ CM, then $\Delta f=f$. Furthermore, $f$ must
be of the following
form $f(z) = 2^zh(z)$, where $h$ is a periodic entire function with period 1.}\\

In 2013, under the restriction on the order of meromorphic function,
Z.X. Chen and H.X. Yi \cite{CY} deduced a uniqueness theorem of
meromorphic functions share three distinct values with their difference operators as follows.\\

{\bf Theorem B}. {\it Let $f$ be a transcendental meromorphic
function such that its order of growth $\rho(f)$  is not an integer
or infinite, let $c\in \mathbb{C}$. If $f$ and $\Delta f (\not\equiv0)$ share three distinct values $e_1,~e_2$, $\infty$, then $f(z+c)=2f(z)$.}\\

In the same paper, Z.X. Chen and H.X. Yi conjectured that the
restriction on the order of growth of $f$ in Theorem B can be
omitted. Clearly, Theorem A showed that the conjecture is right if
$f$ is an entire function of finite order. In the present paper, we
still focus on the conjecture and prove that it holds if
$f$ is a meromorphic function of finite order. In fact, our result is stated as follows.\\

\textbf{Main theorem.} Let $f$ be a transcendental meromorphic
function of finite order, let $\triangle f=f(z+c)-f(z)(\not\equiv0)$, where $c$
is a finite number. If $\triangle f$ and $f$ share three distinct
values $e_1,~ e_2$, $\infty$ CM, then $f=\triangle f$.\\

\textbf{Remark 1.} We point out that there exist meromorphic functions satisfying the conditions of Main theorem. For example,
$f(z)=e^{zln 2}tan(\pi z)$. Obviously, $f=\triangle f=f(z+1)-f(z)$. So, $f$ and $\triangle f$ share $e_1,~e_2$ and $\infty$ CM.\\

\textbf{Remark 2}. The number of shared values cannot be reduced to two. For example, $f(z)=e^{\pi iz}$ and $\triangle f=f(z+1)-f(z)=-e^{\pi iz}$ share $0,~\infty$ CM.
But $f\neq \triangle f$. The example can be seen in \cite{ZL}.\\

\textbf{Remark 3}. Obviously, our main theorem is an improvement of
Theorem A. We also remark that our proof is based on Borel's lemma
\cite{G}. We assume that the reader is familiar with the standard
notations in the Nevanlinna theory, see (\cite{YY, YL}).

\section{Some lemmas}
To prove our result, we recall the difference analogue of the second
main theorem in the value distribution theory.
\begin{lemma}\cite[Theorem 2.4]{HK0} \label{le1}
Let $c\in \mathbb{C}$, let $f$ be a meromorphic function of finite
order with $\Delta f\neq 0$. Let $q\geq 2$, and let $a_1$, $\cdots$,
$a_q$ $\in S(f)$ be distinct periodic functions with period $c$.
Then
$$
m(r,f)+\sum_{i=1}^q m(r,\frac{1}{f-a_i})\leq
2T(r,f)-N_{pair}(r,f)+S(r,f),
$$
where $N_{pair}(r,f)=2N(r,f)-N(r,\Delta f)+N(r\frac{1}{\Delta f})$,
and the exceptional set associated with $S(r, f )$ is of finite
logarithmic measure.
\end{lemma}

A version of Borel's lemma is also needed.
\begin{lemma}\cite[p. 69-70]{G} \label{le2}
Suppose that $n\geq 2$ and let $f_1,~f_2, \cdots, f_n$ be
meromorphic functions and $g_1,~g_2, \cdots, g_n$ be entire
functions such that

(1) $\sum_{j=1}^nf_je^{g_j}=0$,

(2)  when $1\leq j<k\leq n$, $g_j-g_k\neq 0$,

(3) $T(r, f_j)=o(T(r, \exp\{g_h-g_k\}))~(r\rightarrow
\infty,~r\not\in E)$, $E\subset [1,~+\infty)$ of finite logarithmic
measure. Then $f_j=0$ for all $j\in\{1,~2,~\cdots,~n\}$.
\end{lemma}

\section{Proof of Main theorem}
Note that $f,~\triangle f$ share $e_1, ~e_2,$ $\infty$ CM and $f$
is of finite order. Then, there exist two polynomials $\alpha, ~\beta$
such that
\begin{equation}\label{2.1}
\frac{f-e_1}{\triangle f-e_1}=e^{\alpha},~~\frac{f-e_2}{\triangle
f-e_2}=e^{\beta}.
\end{equation}
If $e^{\alpha}=1$ or $e^{\beta}=1$, then $f=\triangle f$. If
$e^{\alpha}=e^{\beta}$, then
$$\frac{f-e_1}{\triangle f-e_1}=\frac{f-e_2}{\triangle f-e_2},
$$
which implies that $f=\triangle f$.

On the contrary, suppose that $f\neq \triangle f$. Then
\begin{equation}\label{2.2}
e^{\alpha}\neq 1,~~e^{\beta}\neq 1,~~e^{\alpha}\neq e^{\beta}.
\end{equation}
Our aim below is to derive a contradiction.

We derive the following expressions from (\ref{2.1}):
\begin{equation}\label{2.3}
f=e_1+(e_2-e_1)\frac{e^{\beta}-1}{e^{\gamma}-1},~~\triangle
f=e_2+(e_2-e_1)\frac{1-e^{-\alpha}}{e^{\gamma}-1},
\end{equation}
where $\gamma=\beta-\alpha$.\\

Obviously,
\begin{equation}\label{2.4}
T(r,f)\leq T(r, e^{\beta})+T(r, e^{\gamma})+S(r,f).
\end{equation}

By the form of $\triangle f$, we have
\begin{equation}\label{2.5}
\begin{aligned}
\triangle f=e_2+(e_2-e_1)\frac{1-e^{\gamma-\beta}}{e^{\gamma}-1}&=(e_2-e_1)(\frac{e^{\beta(z+c)}-1}{e^{\gamma(z+c)}-1}-\frac{e^{\beta}-1}{e^{\gamma}-1})\\
&=(e_2-e_1)(\frac{\beta_1e^{\beta}-1}{\gamma_1e^{\gamma}-1}-\frac{e^{\beta}-1}{e^{\gamma}-1}),
\end{aligned}
\end{equation}
where $\beta_1(z)=e^{\beta(z+c)-\beta(z)}$ and
$\gamma_1(z)=e^{\gamma(z+c)-\gamma(z)}$ are small functions of
$e^\beta$ and
$e^\gamma$, respectively.\\

We claim that $\deg \beta =\deg \gamma$.\\

If $\deg \beta<\deg \gamma$, then $e^\beta$ is a small function of
$e^\gamma$. Suppose that $z_0$ is a zero of $\gamma_1e^{\gamma}-1$,
not a zero of $\beta_1e^{\beta}-1$. Then, it follows from
(\ref{2.5}) that $z_0$ is also a zero of $e^\gamma-1$. Then $z_0$ is
a zero of $\gamma_1-1$. If $\gamma_1-1\neq0$, then
$$
\begin{aligned}
T(r, e^{\gamma})&=\overline{N}(r,
\frac{1}{\gamma_1e^{\gamma(z)}-1})+S(r,e^{\gamma})\\
&\leq N(r, \frac{1}{\beta_1e^{\beta(z)}-1})+N(r,
\frac{1}{\gamma_1-1})+S(r,e^{\gamma})=S(r,e^{\gamma}),
\end{aligned}
$$
a contradiction. Thus, $\gamma_1(z)=e^{\gamma(z+c)-\gamma(z)}=1$. It
means that $\deg \gamma=1$. Note that $\deg \beta<\deg \gamma$, so
$\beta$ is a constant, and say $A$. Thus, again by (\ref{2.5}), we
derive that
$$
\triangle
f=(e_2-e_1)(\frac{\beta_1e^{\beta}-1}{\gamma_1e^{\gamma}-1}-\frac{e^{\beta}-1}{e^{\gamma}-1})=(e_2-e_1)\frac{A-A}{e^\gamma-1}=0,
$$
a contradiction.\\

If $\deg \beta>\deg \gamma$, then $e^{\gamma}$ is a small function
of $e^{\beta}$. Assume that $z_0$ is zero of $e^\beta-1$, not a zero
of $e^\gamma-1$. Then, $z_0$ is a zero of $f-e_1$. Note that $f$ and
$\triangle f$ share $e_1$ CM, so $z_0$ is also a zero of $\triangle
f -e_1$. Put $z_0$ into last form of $\triangle f$ in
(\ref{2.5}), we have
$$
e_1=(e_2-e_1)\frac{\beta_1-1}{\gamma_1e^{\gamma}-1}\big|_{z_0}.
$$
Obviously, $e_1=(e_2-e_1)\frac{\beta_1-1}{\gamma_1e^{\gamma}-1}$.
Otherwise,
$$
\begin{aligned}
T(r, e^{\beta})&=\overline{N}(r,
\frac{1}{e^{\beta}-1})+S(r,e^{\beta})\\
&\leq N(r, \frac{1}{e^{\gamma}-1})+N(r,
\frac{1}{(e_2-e_1)\frac{\beta_1-1}{\gamma_1e^{\gamma}-1}-e_1})+S(r,e^{\beta})=S(r,e^{\beta}),
\end{aligned}
$$
a contradiction. Thus,
$$
e_1=(e_2-e_1)\frac{\beta_1-1}{\gamma_1e^{\gamma}-1}.
$$
Rewrite it as
\begin{equation}\label{2.6}
(e_2-e_1)e^{\beta(z+c)-\beta(z)}-(e_2-e_1)=e_1e^{\gamma(z+c)}-e_1.
\end{equation}
We will prove that $\gamma$ is a constant. On the contrary, suppose that $\deg \gamma\geq 1$. Then,
combining (\ref{2.6}) and $\deg \beta>\deg \gamma$, we obtain that
\begin{equation}\label{2.7}
(e_2-e_1)\beta_1=(e_2-e_1)e^{\beta(z+c)-\beta(z)}=e_1e^{\gamma(z+c)},~~e_2-e_1=e_1.
\end{equation}
It implies  $\beta_1=e^{\gamma(z+c)}$. Rewrite (\ref{2.5}) as
$$
\begin{aligned}
&e_2e^\beta(\gamma_1e^\gamma-1)(e^\gamma-1)+(e_2-e_1)(e^\beta-e^\gamma)(\gamma_1e^\gamma-1)\\
=&
(e_2-e_1)[(\beta_1e^\beta-1)e^\beta(e^\gamma-1)-(e^\beta-1)e^\beta(\gamma_1e^\gamma-1)].
\end{aligned}
$$
Rewrite it as
$$
a_0e^{2\beta}+a_1e^\beta+a_2=0,
$$
where $a_0=(e_2-e_1)[\beta_1(e^\gamma-1)-(\gamma_1e^\gamma-1)]$,
$a_1$, $a_2$ are small functions of $e^\beta$. It indicates that
$a_0=0$. Thus,
\begin{equation}\label{2.8}
\beta_1(e^\gamma-1)=\gamma_1e^\gamma-1.
\end{equation}
Put $\beta_1=e^{\gamma(z+c)}$ into (\ref{2.8}), we have
$$
e^{\gamma(z+c)+\gamma(z)}-2e^{\gamma(z+c)}+1=0,
$$
which implies that $\gamma$ is a constant, a contradiction. Thus, we obtain that $\gamma$ is a constant. The form of $f$ shows
that $f$ is an entire function. Then, it follows from Theorem A that
$f=\triangle f$, a contradiction.

Thus, we prove that $$\deg \beta=\deg \gamma\geq 1.$$

Still set $e^{\beta(z+c)}=\beta_1e^\beta$ and
$e^{\gamma(z+c)}=\gamma_1e^{\gamma}$, where $\beta_1,~\gamma_1$ are two
small functions of $e^\beta$ and $e^\gamma$. Then, due to the forms
of $f,~~\triangle f$, a routine calculation leads to
$$
\begin{aligned}
b_0e^{2\gamma}+b_1e^{\beta+2\gamma}+b_2e^{\beta+\gamma}+b_3e^{2\beta}+b_4e^{2\beta+\gamma}+b_5e^{\beta}+b_6e^{\gamma}=0,
\end{aligned}
$$
where
$$
\left\{
\begin{aligned}
b_0(z)& =(e_2-e_1)\gamma_1,~~b_1(z)= -e_2\gamma_1,\\
b_2(z)&=-(e_2-e_1)+e_2(\gamma_1+1)\\
b_3(z) & = (e_2-e_1)(1-\beta_1),b_4(z)=(e_2-e_2)(\beta_1-\gamma_1),\\
b_5(z) & = -e_2,~~b_6(z)=-(e_2-e_1).\\
\end{aligned}
\right.
$$
Obviously, $b_i$ $(0\leq i\leq 6)$ are small functions of $e^\beta$
and $e^\gamma$. (In fact, for the proof of this result, we just need
the specific forms of $b_0$ or $b_6$.) Rewrite it as
$$
\sum_{i=0}^6 b_ie^{g_i}=0,
$$
where
$$
\left\{
\begin{aligned}
g_0(z) & = 2\gamma,\\
g_1(z) & = \beta+2\gamma, ~~g_2(z)=\beta+\gamma,\\
g_3(z) & = 2\beta,~~ g_4(z)= 2\beta+\gamma,\\
g_5(z) & =  \beta,~~g_6(z)=\gamma.\\
\end{aligned}
\right.
$$
Suppose that $$\deg (\beta)=\deg (\gamma)=n.$$

We claim for any $0\leq j<i\leq 6$
$$
\deg (g_i-g_j)=n.
$$

In the following, we consider several cases to prove the above
claim.\\

\textbf{Case 1.} $i=6$.\\

It is easy to check that
$$
\begin{aligned}
&\deg(g_6-g_4)=\deg (-2\beta)=n,\\
&\deg(g_6-g_2)=\deg (-\beta)=n,~~\deg(g_6-g_0)=\deg (\gamma)=n.
\end{aligned}
$$
Suppose that $\deg(g_6-g_5)=\deg (\gamma-\beta)<n$. Then
$e^{\gamma-\beta}$ is a small function of $e^\beta$ and $e^\gamma$.
We denote by $N_E(r)$ the counting function of those common zeros of
$e^\beta-1$ and $e^\gamma-1$. We firstly prove that $N_E(r)=
S(r,e^\gamma)$. Otherwise, suppose that $N_E(r)\neq S(r,e^\gamma)$.
Assume that $z_0$ is a common zero of $e^\beta-1$ and $e^\gamma-1$.
Then $z_0$ is a zero of $e^{\gamma-\beta}-1$. If
$e^{\gamma-\beta}-1\neq 0$, then
$$
N_{E}(r)\leq N(r, \frac{1}{e^{\gamma-\beta}-1})=S(r, e^\gamma),
$$
a contradiction. Thus, $e^{\gamma-\beta}-1=0$. It means that
$e^\beta=e^\gamma$. So, the form of $f$ yields that $f$ is a
constant, which is impossible.
Thus, we prove that $N_E(r)=S(r,e^\gamma)$.\\

Without loss of generality, assume that $z_0$ is a zero of
$\gamma_1e^\gamma-1$, not a zero of $\beta_1e^\beta-1$. It follows
from (\ref{2.5}) that $z_0$ is also a zero of $e^\gamma -1$. Then
$z_0$ is a zero of $\gamma_1-1$. If $\gamma_1-1\neq 0$, then
$$
\begin{aligned}
T(r, e^{\gamma})&=\overline{N}(r,
\frac{1}{\gamma_1e^{\gamma}-1})+S(r,
e^{\gamma})\\
&\leq N_{E}(r)+N(r, \frac{1}{\gamma_1-1})+S(r,
e^{\gamma})\\
&\leq T(r, \gamma_1-1)+S(r, e^{\gamma})=S(r, e^{\gamma}),
\end{aligned}
$$
a contradiction. Thus, $\gamma_1=e^{\gamma(z+c)-\gamma(z)}=1$, which
implies that $e^{\gamma(z+c)}=e^{\gamma(z)}$ and $\deg \gamma =1$.
Then $\deg(\beta-\gamma)<1$. It means that $\beta-\gamma$ is a
constant, say $A$. Thus, it follows from
$e^{\gamma(z+c)}=e^{\gamma(z)}$ that
$$
e^{\beta(z+c)-\beta(z)}=e^{\beta(z+c)-\gamma(z+c)-(\beta(z)-\gamma(z))}=e^{A-A}=1,
$$
So $e^{\beta(z+c)}=e^{\beta(z)}$. Then, by (\ref{2.5}) we get
$\triangle f=0$, a contradiction. Thus,
$$\deg(g_6-g_5)=n.$$

Suppose that $\deg(g_6-g_3)=\deg (\gamma-2\beta)<n$. The notation
$N_E(r)$ is defined as  above discussion. We firstly prove that
$N_E(r)=S(r,e^\gamma)$. Otherwise, suppose that $N_E(r)\neq S(r,
e^\gamma)$. Without loss of generality, assume that $z_0$ is a
common zero of $e^\beta-1$ and $e^\gamma-1$. Then
$e^{\gamma(z_0)}=1$ and $e^{\beta(z_0)}=1$. Furthermore,
$e^{\gamma(z_0)-2\beta(z_0)}=1$. Clearly, $e^{\gamma-2\beta}$ is a
small function of $e^\gamma$. If $e^{\gamma-2\beta}-1\neq 0$, then
$$
\begin{aligned}
N_E(r)\leq N(r, \frac{1}{e^{\gamma-2\beta}-1})=S(r, e^{\gamma}),
\end{aligned}
$$
a contradiction. Thus, $e^{\gamma-2\beta}=1$ and
$e^\gamma=e^{2\beta}$. Then,
$$
\begin{aligned}
\triangle
f=e_2+(e_2-e_1)\frac{1-e^{\gamma-\beta}}{e^{\gamma}-1}&=e_2+(e_2-e_1)\frac{e^{\beta}(e^{-\beta}-e^{\gamma-2\beta})}{e^{\gamma}-1}
\\
&=e_2+(e_2-e_1)\frac{1-e^{\beta}}{e^{\gamma}-1}.
\end{aligned}
$$
From the forms of $f$ and $\Delta f$, we have
\begin{equation}\label{2.9}
f-e_1=-(\triangle f-e_2).
\end{equation} Since $f$ and $\triangle f$
share $e_1$ and $e_2$ CM, it follows from (\ref{2.9}) that $e_1,
~e_2$ are two Picard values of $f$. Then, by the second main theorem
(see Lemma \ref{le1}), we get
$$
\begin{aligned}
T(r,f)&\leq
N(r,f)+N(r,\frac{1}{f-e_1})+N(r,\frac{1}{f-e_2})\\
&-2N(r,f)+N(r,\Delta
f)-N(r,\frac{1}{\Delta f})+S(r,f)\\
&\leq N(r,\frac{1}{f-e_1})+N(r,\frac{1}{f-e_2})+S(r,f)\leq S(r,f),\\
\end{aligned}
$$
a contradiction. Thus, $N_E(r)= S(r, e^\gamma)$.\\

Similarly as above, we can deduce that $\gamma_1=1$ and $\deg
\gamma=1$. Then, by $\deg (\gamma-2\beta)<n$ and $\deg \gamma= \deg
\beta$, we can set $e^\beta=AH$, $e^\gamma=BH^2$ and
$e^{\gamma-\beta}=CH$, where $A,~B,~C$ are three nonzero constants.
From (\ref{2.5}), a careful calculation leads to
$$
e_2e^\gamma-(e_2-e_1)e^{\gamma-\beta}-e_1=(e_2-e_1)(\beta_1-1)e^\beta.
$$
Rewrite the above equation as
$$
e_2BH^2+b_1H-e_1=0,
$$
where $b_1=-(e_2-e_1)[C+A(\beta_1-1)]$ is a small function of $H$.
Then $e_2=0$ and $e_1=0$. It is impossible. Thus,
$$\deg(g_6-g_3)=\deg (\gamma-2\beta)=n.
$$

Suppose that $\deg(g_6-g_1)=\deg[-(\gamma+\beta)]<n$. The notation
$N_E(r)$ is defined as above discussion. We firstly prove that
$N_E(r)=S(r,e^\gamma)$. Otherwise, suppose that $N_E(r)\neq S(r,
e^\gamma)$. Without loss of generality, assume that $z_0$ is a
common zero of $e^\beta-1$ and $e^\gamma-1$. Then
$e^{\gamma(z_0)}=1$ and $e^{\beta(z_0)}=1$. Furthermore,
$e^{\gamma(z_0)+\beta(z_0)}=1$. Clearly, $e^{\gamma+\beta}$ is a
small function of $e^\gamma$. If $e^{\gamma+\beta}-1\neq 0$, then
$$
\begin{aligned}
N_E(r)\leq N(r, \frac{1}{e^{\gamma+\beta}-1}) =S(r, e^{\gamma}),
\end{aligned}
$$
a contradiction. Thus, $e^{\gamma+\beta}=1$ and
$e^{-\gamma}=e^{\beta}$. Then,
$$
\begin{aligned}
f=e_1+(e_2-e_1)\frac{e^\beta-1}{e^\gamma-1}&=e_1+(e_2-e_1)\frac{e^{-\gamma}-1}{e^\gamma-1}\\
&=e_1+(e_1-e_2)e^{-\gamma},
\end{aligned}
$$
and
$$
\begin{aligned}
\triangle
f=e_2+(e_2-e_1)\frac{1-e^{\gamma-\beta}}{e^{\gamma}-1}&=e_2+(e_2-e_1)\frac{1-e^{2\gamma}}{e^{\gamma}-1}
\\
&=e_1+(e_1-e_2)e^{\gamma}.
\end{aligned}
$$
Furthermore, by $\triangle f=f(z+c)-f(z)$, we have
$$
e_1+(e_1-e_2)e^{\gamma}=(e_2-e_1)(e^{-\gamma(z+c)}-e^{-\gamma}).
$$
Rewrite it as
$$
e_1+(e_1-e_2)e^{\gamma}=(e_2-e_1)(\gamma_2e^{-\gamma}-e^{-\gamma}),
$$
where $\gamma_2$ is a small function of $e^\gamma$ and
$e^{-\gamma}$. Then, it implies that $e_1-e_2=0$, a contradiction. Thus, $\deg(g_6-g_1)=n$.\\

\textbf{Case 2.} $i=5$.\\

It is obvious from the above discussion that
$$
\begin{aligned}
&\deg(g_5-g_4)=\deg (\beta+\gamma)=n,~~\deg(g_5-g_3)=\deg
(-\beta)=n,\\
&\deg(g_5-g_2)=\deg (-\gamma)=n,~~\deg(g_5-g_1)=\deg (-2\gamma)=n.
\end{aligned}
$$
Suppose that $\deg(g_5-g_0)=\deg (\beta-2\gamma)<n$. The notation
$N_E(r)$ is defined as above discussion. We firstly prove that
$N_E(r)=S(r,e^\gamma)$. Otherwise, suppose that $N_E(r)\neq S(r,
e^\gamma)$. Without loss of generality, assume that $z_0$ is a
common zero of $e^\beta-1$ and $e^\gamma-1$. Then
$e^{\gamma(z_0)}=1$ and $e^{\beta(z_0)}=1$. Furthermore,
$e^{\beta(z_0)-2\gamma(z_0)}=1$. Clearly, $e^{\beta-2\gamma}$ is a
small function of $e^\gamma$. If $e^{\beta-2\gamma}-1\neq 0$, then
$$
\begin{aligned}
N_E(r)\leq N(r, \frac{1}{e^{\beta-2\gamma}-1})=S(r, e^{\gamma}),
\end{aligned}
$$
a contradiction. Thus, $e^{\beta-2\gamma}=1$ and
$e^\beta=e^{2\gamma}$. Then,
$$
\begin{aligned}
f=e_1+(e_2-e_1)\frac{e^\beta-1}{e^\gamma-1}&=e_1+(e_2-e_1)\frac{e^{2\gamma}-1}{e^\gamma-1}\\
&=e_2+(e_2-e_1)e^{\gamma},
\end{aligned}
$$
and
$$
\begin{aligned}
\triangle
f=e_2+(e_2-e_1)\frac{1-e^{\gamma-\beta}}{e^{\gamma}-1}&=e_2+(e_2-e_1)\frac{1-e^{-\gamma}}{e^{\gamma}-1}
\\
&=e_1+(e_1-e_2)e^{-\gamma}.
\end{aligned}
$$
Furthermore, by $\triangle f=f(z+c)-f(z)$, we have
$$
e_1+(e_1-e_2)e^{-\gamma}=(e_2-e_1)(e^{\gamma(z+c)}-e^{\gamma}),
$$
which implies that $e_1-e_2=0$, a contradiction. Thus,
$$
\deg (g_5-g_0)=\deg
(\beta-2\gamma)=n.
$$

\textbf{Case 3.} $i=4$.\\

It is obvious from the above discussion that
$$
\begin{aligned}
&\deg(g_4-g_3)=\deg (\gamma)=n,~~\deg(g_4-g_2)=\deg
(\beta)=n,\\
&\deg(g_4-g_1)=\deg (\beta-\gamma)=n, ~~\deg(g_4-g_0)=\deg
(2\beta-\gamma)=n.
\end{aligned}
$$

\textbf{Case 4.} $i=3$.\\

It is obvious from the above discussion that
$$
\begin{aligned}
&\deg(g_3-g_2)=\deg (\beta-\gamma)=n,\\
&\deg (g_3-g_1)=\deg (\beta-2\gamma)=n,~~\deg (g_3-g_0)=\deg
2(\beta-\gamma)=n.
\end{aligned}
$$

\textbf{Case 5.} $i=2$.\\

It is obvious from the above discussion that
$$
\deg(g_2-g_1)=\deg (\gamma)=n,~~ \deg(g_2-g_0)=\deg
(\beta-\gamma)=n.
$$

\textbf{Case 6.} Obviously, $\deg(g_1-g_0)=\deg (\beta)=n$.\\

Thus, the claim is proved. Then, by a Borel' lemma (see Lemma \ref{le2}), we get $b_j=0$
for $0\leq j\leq 6$. But $b_6=-(e_2-e_1)\neq 0$, a contradiction.\\

Thus, the proof of this theorem is finished.

\end{document}